\documentclass[a4paper,12pt]{article}
\date{}

\usepackage{amsfonts}
\usepackage{amsthm}
\usepackage{amsmath}
\usepackage{amsfonts}
\usepackage{latexsym}
\usepackage{amssymb}
\newtheorem{fed}{Definition}[section]
\newtheorem{teo}[fed]{Theorem}

\newtheorem{lem}[fed]{Lemma}
\newtheorem{cor}[fed]{Corollary}
\newtheorem{pro}[fed]{Proposition}
\theoremstyle{definition}
\newtheorem{rem}[fed]{Remark}
\newtheorem{rems}[fed]{Remarks}

\def\noi{\noindent}
\def\QED{\hfill $\blacksquare$}

\def\eps{\varepsilon}

\def\la{\lambda}

\def\cB{\mathcal{B}}

\def\ede{\mathcal{D}}
\def\edemas{\mathcal{D}^+}
\def\edep{\mathcal{P}(\mathcal{D})}
\def\edepfin{\mathcal{P}_0(\mathcal{D})}
\def\edepfins{\mathcal{P}_{0, \ese}(\mathcal{D})}
\def\edepn{\mathcal{P}(\mathcal{D}_n)}
\def\edepm{\mathcal{P}(\mathcal{D}_m)}

\def\cH{\mathcal{H}}

\def\cP{\mathcal{P}}
\def\cQ{\mathcal{Q}}
\def\ere{\mathcal{R}}
\def\ese{\mathcal{S}}
\def\ete{\mathcal{T}}
\def\eme{\mathcal{M}}
\def\ene{\mathcal{N}}

\def\bp{\veebar}

\def\orto{^\perp}
\def\inc{\subseteq}
\def\sii{ if and only if }
\def\inv{^{-1}}
\def\*A{\#\sb A}

\def\eps{\varepsilon}

\def\H{{\cal H}}

\def\glh{Gl(\cH)}
\def\ca{L(\H ) }
\def\cam{L(\H )^+ }

\def\cH{{\cal H}}
\def\dem{\proof  \ }

\def\cP{{\cal P}}

\def\cQ{{\cal Q}}
\def\cS{{\cal S}}
\def\cT{{\cal T}}
\def\cM{{\cal M}}
\def\cN{{\cal N}}

\def\PA{P_{A, \cS}}
\def\PD{P_{D, \cS}}

\def\PAS{\cP(D, \cS)}

\def\rai{^{1/2}}
\def\mrai{^{-1/2}}

\def\api{\langle}
\def\cpi{\rangle}

\def\noi{\noindent}

\def\bm{\left(\begin{array}}
\def\em{\end{array}\right)}
\def\ben{\begin{enumerate}}
\def\een{\end{enumerate}}
\def\barr{\begin{array}}
\def\earr{\end{array}}
\def\iiff{if and only if }
\def\inv{^{-1}}

\def\H{{\cal H}}
\def\lh{{L(\H)}}
\def\lh+{{\lh^+}}

\def\la{\lambda}
\def\eps{\varepsilon}

\def\com{$(A, \cS)$ is compatible}
\def\comd{$(D, \cS)$ is compatible}

\DeclareMathOperator{\Preal}{Re} \DeclareMathOperator{\Pim}{Im}

\DeclareMathOperator*{\convsotipre}{\nearrow}

\DeclareMathOperator{\cc}{\mbox{co}}

\newcommand{\pint}[1]{\displaystyle \left \langle #1 \right\rangle}

\newcommand{\hil}{\mathcal{H}}
\newcommand{\op}{L(\mathcal{H})}

\newcommand{\posop}{L(\mathcal{H})^+}

\newcommand{\bon}{{\mathcal B} = \{e_k\}_{k\in \mathbb {N}}}

\newcommand{\bol}[1]{(#1)_1}
\newcommand{\dist}[2]{\mbox{d}\left(#1,\, #2\right)}
\newcommand{\spec}[1]{\sigma\left(#1\right)}
\newcommand{\pas}[2]{P_{#1,\,#2}}
\newcommand{\Pas}[2]{\mathcal{P}(#1,\,#2)}
\newcommand{\angd}[2]{c_0\left[\,#1,\,#2\,\right]}
\newcommand{\angf}[2]{c\left[\,#1,\,#2\,\right]}
\newcommand{\sang}[2]{s\left[\,#1,\,#2\,\right]}
\newcommand{\con}[2]{K\left[\,#1,\,#2\,\right]}

\newcommand{\generado}[1]{\mbox{span}\left\{#1\right\}}
\newcommand{\conv}{\xrightarrow[n\rightarrow\infty]{}}

\newcommand{\convsot}{\xrightarrow[n\rightarrow\infty]{\mbox{\tiny{S.O.T.}}}}

\newcommand{\convsoti}{\convsotipre_{n\rightarrow\infty}^{\mbox{\tiny
\textbf SOT}}}

\newcommand{\convwot}{\xrightarrow[n\rightarrow\infty]{\mbox{\tiny{W.O.T.}}}}

\newcommand{\chinchulin}[2]{\overline{\chi}_{#1,\,#2}}

\begin{document}


\title{Weighted projections and Riesz frames{}}
\author {Jorge Antezana \thanks{Depto. de Matem\'atica, FCE-UNLP,  La Plata, Argentina
and IAM-CONICET  (e-mail: antezana@mate.unlp.edu.ar)} \and
Gustavo Corach \thanks {Depto. de Matem\'atica, FI-UBA, Buenos Aires, Argentina and IAM-CONICET
(e-mail: gcorach@fi.uba.ar). Partially supported by
CONICET (PIP 2083/00), Universidad de Buenos Aires (UBACYT X050) and
ANPCYT (PICT03-09521)} \and Mariano Ruiz \thanks{Depto. de Matem\'atica, FCE-UNLP,  La Plata, Argentina
and IAM-CONICET  (e-mail: maruiz78@yahoo.com.ar)}
\and Demetrio Stojanoff
\thanks {Depto. de Matem\'atica, FCE-UNLP, 1 y 50 (1900), La Plata, Argentina  and IAM-CONICET
(e-mail: demetrio@mate.unlp.edu.ar). Partially supported CONICET (PIP 2083/00),
Universidad de La PLata (UNLP 11 X350) and ANPCYT (PICT03-09521)}}
\maketitle

\begin{abstract}
Let $\mathcal{H}$ be a (separable) Hilbert space and
$\{e_k\}_{k\geq 1}$ a fixed orthonormal basis of $\mathcal{H}$.
Motivated by many papers on scaled projections, angles of
subspaces and oblique projections, we define and study the notion
of compatibility between a subspace and the abelian algebra of
diagonal operators in the given basis. This is used to refine previous
work on scaled projections, and to obtain
a new characterization of Riesz frames.
\end{abstract}

\bigskip
\noi{\bf Keywords:}
scaled projection, weighted projection, compatibility, angles, frames, Riesz frames.

\noi{\bf AMS classification:}
Primary 47A30, 47B15
\bigskip

\section{Introduction}

Weighted projections (also called scaled projections)
play a relevant role in a
variety of least-square problems. As a sample of their
applications and of their relatives, namely, weighted pseudoinverses, 
they have been
used in optimization (feasibility theory,
interior point methods), statistics (linear
regression, weighted estimation), and signal processing
(noise reduction).

Frequently, \  weighted \ pseudoinverses \ take \ the forms
$(ADA^t)\inv AD$, $(ADA^t)^\dagger AD$, $(ADA^*)\inv AD$ or $(ADA^*)^\dagger AD$,
according to  the field (real or complex) which is involved in the problem and to 
different hypothesis of invertibility. 
Analogous formulas hold
for the corresponding weighted projections. In general $D$ is
a positive definite matrix and $A$ is a full column rank matrix.

In a series of papers, Stewart \cite{[stewart]}, O'Leary
\cite{[O'L]}, Ben-Tal and Taboulle \cite{[bental]}, Hanke and
Neumann \cite{[HaNe]}, Forsgren \cite{[Fgren1]}, Gonzaga and Lara
\cite{[clovislara]}, Forsgren and Sp\"orre \cite{[Fgren2]}, and Wei
\cite{[wei3]}, \cite{[wei1]}, \cite{[wei2]}, \cite{[wei4]} have
studied and computed quantities of the type
$$
\sup _{D\in \Gamma } \|\gamma (D, A)\|,
$$
where $\Gamma $ denotes a certain subset of positive definite invertible
matrices and $\gamma (D, A)$ is any of the weighted  pseudoinverses
mentioned above. The reader is referred to the papers by
Forsgren and Sp\"orre \cite{[Fgren1]}, \cite{[Fgren2]} for  excellent surveys on the
history and motivations of the problem of estimating the supremum above.
It should be said,
however, that the references mentioned above only deal with the finite
dimensional context. In order to deal with increasing dimensions or arbitrarily
large data sets, we  present the problem in an infinite dimensional
Hilbert space.

Moreover, we present a different approach to this theory, valid
also in the finite dimensional context, based on techniques and results on
generalized selfadjoint projections. Recall that, if $D$ is a
selfadjoint operator on a complex (finite or infinite dimensional)
Hilbert space $\H$, another operator $C$ on $\H$ is called
$D$-selfadjoint if $C$ is Hermitian with respect to the
Hermitian sesquilinear form
$$
\api \xi , \eta \cpi_D = \api D\xi , \eta \cpi  \quad
(\xi , \eta \in \H),
$$
i.e. if $DC= C^*D$. We say that a closed subespace $\ese$ of $\H$ is
$compatible$ with $D$ (or that the pair $(D, \ese )$ is compatible) if
there exists an $D$-selfadjoint projection $Q$ in $\H$ with image $\ese$.
It is well known \cite{[CMS1]} that in finite dimensional spaces,
every subspace is compatible with any positive semidefinite operator $D$.
In infinite dimensional spaces this is not longer true; however,
every (closed) subspace is compatible with any positive invertible
operator and, in general, compatibility can be characterized in
terms of angles between certain closed subspaces of $\H$, e. g.,
the angle between $\ese$ and $(D\ese )\orto$.

If the pair $(D, \ese )$ is compatible, the set of
$D$-selfadjoint projections onto $\ese$ may be infinite; nonetheless, a distinguished
one denoted by $\PD$, can be defined and computed
(see \cite{[CMS1]} or section 2.2 below).

In the finite dimensional case, from the point of view of $D$-selfadjoint projections, 
the study of weighted projections allows us to 
obtain simpler proofs of some known results. Another advantage 
this perspective offers is
that these proofs can be easily extended to more general settings
which are also important in applications. These applications include projections with
complex weights and the infinite dimensional case.
Moreover, this approach establishes the relationships
among the quantities that have appeared throughout the study 
on weighted projection (usually, operator norms,
vector norms and angles).

A well known result due to Ben-Tal and Teboulle
states that the solutions to
weighted least squares problems lies in the convex hull of
solutions to some non-singular square subsystems. We refer the reader
to Ben-Tal and Teboulle's paper
\cite{[bental]}, or \cite{[Fgren1]}, \cite{[wei3]} for
the following formulation:
let $A$ be an $m\times n$ matrix of full rank. Denote by  $J(A)$ the set of all
$m\times m$ orthogonal diagonal projections   such that
$QA:\mathbb{C}^n\rightarrow R(Q)$ is bijective. Then, for every
$m\times m$ positive diagonal matrix $D$,
\begin{equation}\label{BeTa}
A(A^*DA)^{-1}A^*D= \sum_{Q\in
J(A)}\left(\frac{\det(D_Q)|\det(A_Q)|^2} {\sum_{P\in J(A)}
\det(D_P)|\det(A_P)|^2}\right) A(QA)^{-1}Q
\end{equation}
where $A_Q$ (resp. $D_Q$) is $QA$ (resp. $QD$) considered as a
square submatrix of $A$ (resp. $D$).

In section 3 we show that, if  $\ese = R(A)$, then for every
$D\in \edemas$ and $Q\in J(A)$  the following identities hold:
$$
A(A^*DA)^{-1}A^*D=\pas{D}{\ese} \quad \hbox{ and } \quad A(QA)^{-1}Q=\pas{Q}{\ese},
$$
where $\pas{D}{\ese} $ and $\pas{Q}{\ese}$ denote the distinguished projections onto $\ese$ which are $D$-selfadjoint
and $Q$-selfadjoint, respectively.

Then, Ben-Tal and Teboulle's formula (\ref{BeTa})
can be rewritten in the following way: if $R(A) = \ese$ and for every $D\in \edemas$,
\[
\pas{D}{\ese}\in \cc \{\pas{Q}{\ese}:Q\in J(A)\}.
\]
This  implies, in particular,  that $\displaystyle
\sup_{D\in\edemas_n}\|\pas{D}{\ese}\|\leq \max_{Q\in
J(A)}\|\pas{Q}{\ese}\|.$ The same inequality was proved
independently by O'Leary in \cite{[O'L]}, while the reverse
inequality was initially proved by Stewart \cite{[stewart]}. A slight generalization of
Stewart's result is proved in this section.
Another application of the projections techniques provides an easy proof of a result of
Gonzaga and Lara \cite{[clovislara]} about scaled projections, even for complex weights.

In section 4, we extend the notion of compatibility of a closed subspace, with respect
to certain subsets of $\cam$. Given $\Gamma \inc \cam$ and a closed subspace $\ese$,
we say that $\ese$ is compatible with  $\Gamma$ if $(D, \ese )$ is compatible
for every $D \in \Gamma$
and it satisfies Stewart's condition:
$$\sup _{D\in \Gamma} \| \PD\| < \infty.$$
For a fixed orthonormal basis $\cB = \{e_n\}_{n \in \mathbb {N}}$ of $\H$,
we denote by $\ede$ the diagonal
algebra with respect to $\cB$,  i.e. $D\in \ede$ if $De_n = \la_n e_n$
$(n \in \mathbb {N}$) for a bounded sequence $(\la_n)$ of complex numbers. Next, we consider
compatibility of $\ese$ with respect to
\ben
\item $\edemas$,  the set of positive invertible elements of $\ede$ (i.e. all $\la_n > \eps$,
for some $\eps >0$);
\item $\edep$, the set of projections in $\ede$ (i.e. all $\la_n = 0$ or $1$);
\item $\edepfin$, the set of elements in $\edep$ with finite rank, and
\item $\edepfins$, the set of elements $ Q \in \edepfin$ such that $R(Q) \cap \ese = \{ 0 \}$
\een
  For a  closed subspace $\ese$, we show that  compatibility with any of these sets 
is equivalent. In the first case,  we say that $\ese$ is compatible with 
the basis $\cB$ (or $\cB$-compatible).

This notion is very restrictive. Nevertheless, the class of subspaces which are
compatible with a given basis $\cB$ has its own interest. Indeed, as we show in section 5,
if $\dim \ese\orto = \infty$, then  $\ese$ is $\cB$-compatible \iiff the class of
frames whose preframe operators (in terms of the basis $\cB$)
have nullspace $\ese$, consists of Riesz frames
(see section 5 for definitions or Casazza \cite{[Cas]},
Christensen \cite{[Chris]}, \cite{[Chri]} for
modern treatments of Riesz frame theory and applications).
We completely characterize compatible subspaces with $\cB$ in terms of Friedrichs
angles (see Definition \ref{fried}) and we obtain an
analogue of Stewart-O'Leary identity. 
Let $\ese$ be a closed subspace of $\H$. For $J \inc \mathbb{N}$,
denote by $\H_J $ the closed span of the set $\{ e_n : n \in J\}$ and $P_J$ the orthogonal
projector onto $\H_J$. In the case that $J = \{ 1, \dots , n\}$, we denote $\H_n $ and
$P_n$ instead of $\H_J$ and $P_J$. Then, the main results of this paper are:
\ben
\item The following conditions are equivalent:
\ben
\item $\ese$ is compatible with $\edemas$;
\item $\sup \{ \angf{\ese }{\H_J}: J \inc \mathbb{N}\} < 1 $, were $\angf{\ete}{\eme} $ denotes
the Friedrichs angle between the closed subspaces $\ete$ and $\eme$;
\item $\sup \{ \angf{\ese }{\H_J}: J \inc \mathbb{N} $ and $J$ is finite $\} < 1 $;
\item all pairs $(P_J, \ese)$ are compatible and
$\sup  \{ \|P_{P_J , \ese} \| : J \inc \mathbb{N}\} < \infty$.
\een
In this case
$$
\sup \Big{\{} \| \PD\|: D\in \edemas \Big{\}} = \sup  \Big{\{} \|P_{P_J , \ese} \| : J \inc \mathbb{N}\Big{\}}
= \Big{(}1- \sup_{Q\in\edep}\angf{\ese}{R(Q)}^2 \Big{)}\mrai .
$$
\item $\ese$ is compatible with $\edemas$ \iiff
\ben
\item $\displaystyle\ese = \overline{\cup _{n \in \mathbb{N}} \ \ \ese \cap \H_n}$ and
\item for every $n \in \mathbb{N}$, the subspace  $\ese \cap \H_n$ is compatible with $\cB$ and
there exists $M > 0$ such that  $\sup\{ \|\pas{P_J}{\ese\cap\H_n}\| :  J \inc \mathbb{N}\}
\leq M$ for every $n \in \mathbb{N}$.
\een
\item If $\dim \ese <\infty$, then $\ese$ is compatible with $\cB$ \iiff there exists
$n \in \mathbb{N}$ such that $\ese \inc \H_n$.
\een

\section{Preliminaries}

Let $\hil$ be a separable Hilbert space and $\ca$ be the
algebra of bounded linear operators on $\hil$. For an operator $A
\in \ca$, we denote by $R(A)$ the range or image of $A$, $N( A)$ the
nullspace of $A$, $\sigma (A)$ the spectrum of $A$, $A^*$ the adjoint
of $A$, $\rho(A)$ the spectral radius of $A$, $\|A\|$ the
usual norm of $A$ and, if $R(A)$ is closed, $A^\dagger$ the
Moore-Penrose pseudoinverse of $A$.

Given a closed subspace $ \ese $ of $ \hil$, we denote by $P_\ese$ the
orthogonal (i.e. selfadjoint) projection onto $\ese$. If $B \in
\ca $ satisfies $P_\ese B P_\ese = B$, we consider the compression
of $B$ to $\ese$, (i.e. the restriction of $B$ to $\ese$ as a
map from $\ese$ to $\ese$), and we say that 
$B$ is considered as $acting$ on $\ese$. 

Given a subspace $\ese$ of $\hil$, its unit ball is denoted by
$\bol{\ese}$, and its closure by $\overline{\ese}$.
The distance between two subsets $S_1$ and $S_2$
of $\hil$ is
\[
\dist{S_1}{S_2}=\inf\{\|x-y\|:\;x\in S_1 \;\;y\in S_2\}.
\]
Along this note we use the fact that every subspace $ \ese $ of $
\hil$ induces a representation of elements of $\op$ by $2 \times
2$ block matrices. We shall identify each  $A\in\op$ with
a $2\times 2$-matrix
$
\begin{pmatrix}
  A_{11} & A_{12} \\
  A_{21} & A_{22}
\end{pmatrix}\begin{array}{l}
  \ese  \\
  \ese^\bot
\end{array}
$, which we write to emphasize the decomposition which induces it.  Observe that $
\begin{pmatrix}
  A_{11}^* & A_{21}^* \\
  A_{12}^* & A_{22}^*
\end{pmatrix}$ is the matrix which represents $A^*$.


\subsection{Angle between subspaces}

Among different notions of angle between subspaces in a Hilbert
space, we consider two definitions due to Friedrichs
and Dixmier (see \cite{[Di]} and \cite {[Fr]}).

\bigskip
\begin{fed}[Friedrichs] \label{fried}\rm
Given two closed subspaces $\eme$ and $\ene$, the {\it angle}
between $\eme$ and $\ene$ is the angle in $[0,\pi/2]$ whose cosine
is defined by
\[
\angf{\eme}{\ene}=\sup\{\,|\pint{\xi , \, \eta}|:\;\xi\in
\eme\ominus (\eme\cap \ene), \;\eta\in \ene\ominus (\eme\cap
\ene)\;\mbox{and}\;\|\xi\|=\|\eta\|=1 \} .
\]
Then, the sine of this angle is
$$
\sang{\eme}{\ene} = (1-\angf{\eme}{\ene})\rai = \dist{\ \bol{\eme}}{\ene\ominus(\eme\cap \ene)} .
$$
The last equality follows from the definition by direct computations.
\end{fed}


\bigskip
\begin{fed}[Dixmier]\rm
Given two closed subspaces $\eme$ and $\ene$, the {\it minimal
angle} between $\eme$ and $\ene$ is the angle in $[0,\pi/2]$ whose
cosine is defined by
\[
\angd{\eme}{\ene}=\sup\{\,|\pint{\xi , \, \eta}|:\;\xi\in \eme,
\;\eta\in \ene\;\mbox{and}\;\|\xi\|=\|\eta\|=1 \}
\]
\end{fed}

\noi
The reader is referred to the excellent survey by F. Deutsch \cite {[De]}
and the book of Ben-Israel and Greville \cite {[BG]} which also have
complete references. The next  propositions collect the results about
angles which are relevant to our work.

\bigskip
\begin{pro}\label{propiedades elementales de los angulos}
Let $\eme$ and $\ene$ be to closed subspaces of $\hil$.
Then \ben
  \item $0\leq \angf{\eme}{\ene}\leq \angd{\eme}{\ene}\leq 1$
  \item $\angf{\eme}{\ene}=\angd{\eme\ominus(\eme\cap
\ene)}{\ene}=\angd{\eme}{\ene\ominus(\eme\cap \ene)}$
\item $\angf{\eme}{\ene}= \angf{\eme^\bot}{\ene^\bot}$
  \item $\angd{\eme}{\ene}=\|P_\eme P_\ene\|=\|P_\eme P_\ene
P_\eme\|^{1/2}$
  \item $\angf{\eme}{\ene}=\|P_\eme P_\ene-P_{(\eme\cap
\ene)^\bot}\|=\|P_\eme P_\ene P_{(\eme\cap \ene)^\bot}\|$

\item The following statements are equivalent
\ben
  \item[i.] $\angf{\eme}{\ene}<1$
  \item[ii.] $\eme+\ene$ is closed
  \item[iii.] $\eme^\bot + \ene^\bot$ is closed
\een \een
\end{pro}

\bigskip
\begin{pro}\label{producto con rango cerrado}(\cite {[Bo]}, \cite {[Iz]})
Given $A,B\in\op$, then $R(AB)$ is closed \sii \  $\angf{R(B)}{N(A)}<1$.
\end{pro}

\bigskip
\begin{pro}\label{aproximacion a la interseccion}
Let $P$ and $Q$ be two orthogonal projections defined on $\hil$.
Then,
\[
\|(P Q)^k-P\wedge Q\|=\angf{R(P)}{R(Q)}^{2k-1}
\]
where $P\wedge Q$ is the orthogonal projection onto $R(P)\cap R(Q)$.
\end{pro}

\bigskip
\begin{pro}[Ljance-Ptak \cite{[ptak]}]\label{formula de ptak}
Let $Q$ be a projection with range $\ere$ and with nullspace
$\ene$. Then
\begin{align*}
\|Q\|=\frac{1}{\Big(1-\|P_\ere\,P_\ene\|^2\Big)^{1/2}}=
\frac{1}{\Big(1-\angf{\ere}{\ene}^2\Big)^{1/2}}=\sang{\ere}{\ene}^{-1}
\end{align*}
\end{pro}

\subsection{$D$-selfadjoint projections and compatibility}\label{oblique}

Any selfadjoint operator $D \in \ca $ defines a bounded Hermitian
sesquilinear form
$\api\xi, \eta \cpi _D = \api D \xi , \eta \cpi$, $ \xi , \ \eta \in \H $.
The  $D$-orthogonal subspace of a subset $\cS$ of $\H$  is
$\cS^{\perp _D} :=
\{\xi: \ \api D \xi, \eta \cpi \  = 0 \ \ \forall \eta \in \cS\ \} =
D\inv (\cS ^\perp )=D(\cS)^\perp .
$

We say that $C\in \ca$ is $D$-$selfadjoint$ if $DC = C^*D$.
Consider the set of $D$-selfadjoint projections whose range is exactly $\cS$:
$$
\PAS = \{Q \in \cQ: R(Q) = \cS, DQ = Q^*D\}.
$$
A pair $(D, \cS ) $ is called $compatible$ if $\PAS$ is not empty.
Sometimes we say that $D$ is $\cS$-compatible or that $\cS$ is
$D$-compatible.

\bigskip
\begin{rem}\label{douglas}\rm
It is known (see Douglas \cite{[Do]}), that if $\H_1 , \H_2 , \H_3$ are Hilbert spaces,
$B\in L(\H _3, \H_2 ) $ and $ C \in L(\H _1, \H_2 )$,
then the following conditions are equivalent: \\
a) $R(B) \inc R(C)$; \\
b) there exists a positive number $\lambda $ such that $BB^* \le \lambda CC^*$ and \\
c) there exists $A \in L(\H _1, \H_3 )$ such that $B = CA$.

\noindent
Moreover, there exists a unique operator $A$ which  satisfies
the conditions $ B= CA$  and $ R(A) \inc \overline{R(C^*)}. $
In this case, $N( A) = N(B)$ and $\|A\|^2 = \inf \{\lambda : BB^* \le \lambda CC^* \}$;
$A$ is called the {\it reduced} solution of the equation $CX=B$. If
$R(C)$ is closed, then $A= C^\dagger B$.
\end{rem}

\bigskip
\noindent
In the following theorem we present several results about
compatibility, taken from \cite{[CMS1]} and \cite{[CMS2]}.

\begin{teo} \label{PA}
If $D  \in \ca$ is selfadjoint, and $\cS $ is a closed subspace of $\H$,
we denote  $D= \bm {cc} a & b\\b^* & c \em  \barr{l} \cS \\\cS\orto \earr$. Then:
\ben \item
\comd \ \iiff $R(b) \inc R(a)$ \iiff $\cS + D\inv (\cS\orto ) =  \H$.
\item In this case, if $d\in L(\cS^\perp , \cS)$
is the reduced solution of the equation $ax=b$ then
$$ \PD =  \bm {cc} 1 & d \\ 0 & 0 \em  \in \PAS $$
and, if $\cN = D\inv(S^\perp)\cap \cS$, then
$N(\PD) = D\inv (\cS \orto ) \ominus \cN$.
\item  If $D\in \cam$ then $\cN  = N(D)  \cap \cS$ and,
for every $Q \in \PAS$, there is
$z \in L(\cS^\perp, \cN)$ such that
\begin{equation}\label{1zd}
Q = \PD + z = \bm {ccc} 1 &0& d \\ 0 & 1& z \\ 0&0 & 0
\em \barr{l} \cS\ominus \cN
\\ \cN \\ \cS ^\perp \earr .
\end{equation}
Observe that  $\PAS$ has a unique element (namely, $\PD$) if and only if
$N(D)  \cap \cS = \{ 0 \}$.
\item $\PD$ has minimal norm in $\PAS$, i.e.
$ \|\PD \| = \min \{\ \|Q\| : Q \in \PAS \}.$
\een
\end{teo}

\medskip
\noi
The reader is referred to \cite{[CMS1]},  \cite{[CMS2]} and \cite{[CMS3]} for several applications of
$\PD$ (see also Hassi and Nordstr\"om \cite{[HN]}).

\medskip
\noi
From now on, we shall suppose that $D\in \posop$, in which case
$D^{1/2}$ denotes the positive square root of $D$.

\bigskip
\begin{rem}\label{equiva} \rm
Under additional hypothesis on $D$, other characterizations of
compatibility can be used. We mention a sample of these, taken
from \cite {[CMS1]}  and \cite {[CMS2]}:

\ben
\item
If $R(PDP)$ is closed (or, equivalently, if $R(PD\rai) $ or $D\rai
(\cS )$ are closed), then \comd .
In this case, if $D= \bm {cc} a & b\\b^* & c \em $, then
$\PD = \bm {cc} 1 & a^\dagger b \\ 0 & 0 \em ,
$
since $a = PDP $ has closed range, and $a^\dagger b$ is the
reduced solution of $ax=b$.

\item
If $D $ has closed range then 
 the pair \comd \  $\iff $ \  $R(PDP)$ is closed
\ $\iff $ \  $R(DP)$ is closed $\iff $ $c[N(D)  , \ese] <1$.

\item If $P, Q$ are orthogonal projections with  $R(P) = \ese$,
then $ (Q, \ese )$ is compatible $ \iff$ $R(QP)$ is closed $ \iff $ $ c[N(Q)  , \ese ]< 1 $.
Moreover, if $(Q, \ese ) $ is compatible, then $\H = \ese +  Q \inv (\ese \orto ) = \ese +
(R(Q) \cap \ese \orto ) + N(Q) $ and, if $\cN = N(Q)  \cap \ese
$ and  $\cM = \ese \ominus \cN $, then
$\cM \oplus (N(Q)  \oplus (R(Q) \cap \ese \orto ))=\H$,
and $P_{Q, \ \cM }$ is the projection onto $\cM$ given by this
decomposition.
In particular,  
if $\ese \oplus N(Q)  = \H$, then $P_{Q, \ese } $ is the projection onto
$\ese$ given by this decomposition.
Observe that $P_{Q, \ese } = P_\cN + P_{Q, \ \cM}$. It follows that
$$\|P_{Q,\ese} \| = \|P_{Q,\ \cM} \| = (1-\|(1-Q)P_\cM\|^2)\mrai
= (1-\angf{N(Q)  }{ \ese } ^2)\mrai
= \sang{N(Q)  }{ \ese } \inv.
$$
\een
Observe that, in finite dimensional spaces, every pair $(D,\ese)$ is compatible because
every subspace, a fortiori $R(PDP)$, is closed.
\end{rem}

\bigskip
\noi
We end this section with the following technical result, which we shall
need in what follows:
\bigskip

\begin{pro} \label{compri}
Suppose that $\ese \inc \cT$ are closed subspaces of $\cH$ 
and $P_\cT D = D P_\cT$. If
$R(P_\ese D P_\ese )$ is closed and $D = \bm {cc} D_1 & 0 \\0 &
D_2 \em \barr{l}  \cT \\ \cT\orto \earr $, then $ \PD = \bm {cc}
P_{D_1, \, \ese} & 0 \\ 0 & 0 \em  \barr{l}  \cT \\ \cT\orto \earr
, $ where we consider $P_{D_1, \,  \ese} $ as acting on $\cT$.
\end{pro}
\dem Let $D_1 = \bm {cc} a & b \\ b^* & c \em \barr {l} \ese \\
\cT \ominus \ese \earr$. Then $P_{D_1, \, \ese} = \bm {cc} 1 &
a^\dagger b \\ 0 & 0 \em $. On the other hand, if $$ D = \bm {ccc}
a & b & 0 \\ b^* & c & 0 \\ 0 & 0 & D_2 \em \barr{l} \ese \\ \cT
\ominus \ese \\ \cT \orto \earr , \ \hbox {then } \ \PD = \bm
{ccc} 1 & a^\dagger b & a^\dagger 0 \\ 0 & 0 & 0 \\ 0 & 0 & 0 \em
\barr{l} \ese \\ \cT \ominus \ese \\ \cT\orto \earr = \bm {cc}
P_{D_1, \, \ese} & 0 \\ 0 & 0 \em . $$ \QED

\begin{rem}\rm
Proposition \ref{compri} is still valid if the assumption that
``$R(P_\ese D P_\ese )$ is closed" is replaced by ``the pair \comd ".
The proof follows the same lines
but is a little bit more complicated, because it uses the more general
notion of reduced solutions (see Remark \ref{douglas})
instead of Moore-Penrose pseudoinverses.
One must also show that the pair $(D_1 , \cS)$ is compatible in $L(\cT )$.
\end{rem}

\section{Scaled projections in finite dimensional spaces}

In this section we study scaled projections in finite dimensional
Hilbert spaces from the viewpoint of $D$-selfadjoint projections.
This is a new geometrical approach to the widely studied subject
of weighted projections which may be helpful in the applications.
In the next section we shall use this approach to extend some of
these results to
infinite dimensional spaces. Additionally, we shall
study the projections with \textit{complex weights} as those considered by
M. Wei (\cite {[wei1]}, \cite {[wei2]}, \cite {[wei3]} and
E. Bobrovnikova and S. Vavasis  \cite{[vavasis2]} and we shall
prove generalizations of some well known results about
classical scaled projections to the complex case.

Throughout this section, $\ede_n$ denotes the abelian algebra of
diagonal $n\times n$ complex matrices, $\edemas_n$ denotes the
set of positive invertible matrices of $\ede_n$
and $\edepn$ denotes the set of projections in $\ede_n$.

Scaled projections are connected with scaled pseudo-inverses which
appear
in weighted least squares problems of the form
\[
\min \{\|D^{1/2}(\beta-A\xi)\|^2:\xi\in \mathbb{C}^n\},
\]
where $m \ge n$, $A$ is an $m \times n$ matrix of full rank, $\beta \in  \mathbb{C}^m$ and
$D\in \edemas_m$. 
It is well known that the solution to this problem is
\[
\xi=(A^*DA)^{-1}A^*D \beta .
\]
The operator $A_D^\dagger=(A^*DA)^{-1}A^*D$ is called a
{\it weighted pseudo-inverse \rm} of $A$.

\medskip

In some situations it is useful to have a bound for the norms of
the scaled pseudo-inverses $A_D^\dagger$. In order to study this
problem, G. W. Stewart \cite{[stewart]} used the oblique projections
$P_D=A(A^*DA)^{-1}A^*D$ and proved that
\begin{equation}\label{Stewart en matrices}
M_A = \sup\{\|A(A^*DA)^{-1}A^*D\|:\; D \in \edemas_m \}<\infty.
\end{equation}

For $A\in \mathbb{R}^{n\times n}$ and $I\subseteq \{1,\ldots, n\}$ let $m_I$
denote the minimal non-zero singular value of the submatrix corresponding to the rows
indexed by $I$ of a matrix $U$ whose columns form an orthonormal basis of $N(A)^\bot$.

Stewart \cite {[stewart]} proved that
\begin{equation}\label{StewartOLeary}
M_A \inv  \le \min \{ m_I : I \inc \{ 1, 2, .... , n\}\},
\end{equation}
and O'Leary \cite {[O'L]} proved that both numbers actually coincide.

Independently, A.
Ben Tal and M. Teboulle  \cite{[bental]}
proved the next theorem, which refines Stewart-O'Leary's result:

\bigskip
\begin{teo}\label{bental teboulle}
Let $A$ be an $m\times n$ matrix of full rank and let $J(A)$ be
the set of all  $Q\in \edepm$ such that
$QA:\mathbb{C}^n\rightarrow R(Q)$ is bijective. Then, for every
$D\in \edemas_m$ it holds that
\[
A(A^*DA)^{-1}A^*D=
\sum_{Q\in J(A)}\left(\frac{\det(D_Q)|\det(A_Q)|^2}
{\sum_{P\in J(A)} \det(D_P)|\det(A_P)|^2}\right) A(QA)^{-1}Q
\]
where $A_Q$ (resp. $D_Q$) is $QA$ (resp. $QD$) considered as a
square submatrix of $A$ (resp. $D$). In particular
\[
P_D\in \cc\Big{\{}A(QA)^{-1}Q:\;Q \in J(A) \Big{\}}
\]
\end{teo}

\noi
The reader will find illustrative surveys in the papers by Forsgren  \cite{[Fgren1]}, Forsgren and Sp\"orre
\cite{[Fgren2]} and in Ben-Israel and Greville's book \cite{[BG]} nice surveys on these
matters. See also the papers by Hanke and Neumann (\cite{[HaNe]}), Gonzaga and Lara \cite {[clovislara]},
Wei \cite{[wei4]}, \cite{[wei2]}, \cite{[wei1]} and Wei and de Pierro  \cite{[wei3]}.


\bigskip
\noi
Given a fixed positive diagonal matrix $D
\in \ede_m$, the solution of
\[
\min \{\|D^{1/2}(\beta-\xi)\|^2:\xi\in R(A)\}
\]
is given by $\xi=P_D\, \beta$. Observe that $\|D^{1/2}\cdot\|$ is the norm
induced by the inner product $\pint{D\cdot,\,\cdot}$, and
therefore, $P_D$ is the (unique) projection onto $R(A)$ that is
orthogonal with respect to the inner product $\pint{D\cdot,\,\cdot}$.
Therefore, under the notations of section \ref{oblique},
$P_D=\pas{D}{R(A)}$. It is natural to ask if
$A(QA)^{-1}Q$ coincides with  $\pas{Q}{R(A)}$ for every $Q$ under the
conditions of Ben Tal and Teboulle's Theorem \ref{bental teboulle}. The answer to this
question is the goal of the next Proposition. 

\bigskip
\begin{fed}[ Dixmier \cite{[Di]}]\label{gp}
Two closed subspaces $\ese$ and $\ete$ of $\hil$
are in  position P' if it holds that $\ete \orto \cap \ese=\ete\cap\ese^\bot=\{0\}$.
In this case, we write $\ese \bp \ete$.
\end{fed}

\bigskip
\begin{pro}\label{traduccion}
Given an $m\times n$ matrix  $A$ of full rank, let $D\in
\edemas_m$ and let $Q$ be a diagonal projection. Then
\ben
  \item $\pas{D}{R(A)}=A(A^*DA)^{-1}A^*D$.
  \item  $QA:\mathbb{C}^n\rightarrow R(Q)$ is bijective \iiff $R(Q)\bp R(A)$.
  \item If $R(Q)\bp R(A)$ then $\pas{Q}{R(A)}=A(QA)^{-1}Q$.
\een
\end{pro}

\dem \ben
  \item It suffices to observe  the coincidence of the range (resp. nullspace) of both
  projections  $\pas{D}{R(A)}$ and $A(A^*DA)^{-1}A^*D$.

  \item  $QA:\mathbb{C}^n\rightarrow R(Q)$ is
  bijective if and only if $QA:\mathbb{C}^n\rightarrow R(Q)$ and
  $A^*Q:R(Q)\rightarrow \mathbb{C}^n$ are injective.
  As $QA:\mathbb{C}^n\rightarrow R(Q)$ is
  injective if and only if $R(A)\cap N(Q)=\{0\}$ and
  $A^*Q:R(Q)\rightarrow \mathbb{C}^n$ is injective if and only if
  $R(Q)\cap R(A)^\bot= R(Q)\cap N(A^*)=\{0\}$, it follows that
  $QA:\mathbb{C}^n\rightarrow R(Q)$
  is bijective if and only if $R(Q)\bp R(A)$.

\item Clearly, $A(QA)^{-1}Q$ is a projection whose range is
$R(A)$ and whose nullspace is $N(Q)$. But, $\pas{Q}{R(A)}$ is also
a projector with the same range and nullspace as
  $A(QA)^{-1}Q$. In fact, by item 2 in Theorem \ref{PA}, 
  $$
  R(\pas{Q}{R(A)})=R(A)  \ \ \hbox { and }  \ \
  N(\pas{Q}{R(A)})= Q^{-1}(R(A)^\bot)\ominus(N(Q) \cap R(A))=N( Q).
  $$
  Hence $A(QA)^{-1}Q = \pas{Q}{R(A)}$. \QED
\een

\medskip

\noi Using Proposition \ref{traduccion} we can restate Theorem
\ref{bental teboulle} in the following way:

\bigskip
\begin{teo}\label{bental teboulle con pases}
Let $\ese$ be a subspace of $\mathbb{C}^n$ and let
$D\in\edemas_n$. Then
\begin{equation}\label{bental en matrices}
\pas{D}{\ese}\in \cc \{\pas{Q}{\ese}:Q\in\edepn\;\,
\mbox{and}\;\,R(Q)\bp\ese\}.
\end{equation}
In particular,
\begin{equation}\label{stewart}
\sup_{D\in\edemas_n}\|\pas{D}{\ese}\|\leq
\max\Big{\{}\|\pas{Q}{\ese}\| : Q\in\edepn:\,R(Q)\bp\ese \Big{\}}.
\end{equation}
\end{teo}

\begin{rem}
Inequality (\ref{stewart}) is actually an equality. The converse
inequality was proved by Stewart and it is also a consequence of the next
Proposition which follows essentially Stewart's ideas.
\end{rem}

\bigskip
\begin{pro}\label{pqese menor que pdese}
Let $\ese$ be a subspace of $\mathbb{C}^n$ and denote by
$\edemas_{0,n}$ the set of all diagonal positive semidefinite $ n
\times n$ matrices. Then
\begin{equation}\label{formula con los semidefinidos}
\sup_{D\in\edemas_n}\|\pas{D}{\ese}\|=\sup_{D\in\edemas_{0,n}}\|\pas{D}{\ese}\|.
\end{equation}
\end{pro}

\dem Let $D\in\edemas_{0,n}$  and consider the sequence of
invertible positive operators $\{D_k\}_{k\geq 1}$ defined by
\[
D_k=D+\frac{1}{k} \ I.
\]
If $\ene=\ese\cap N(D)$, then 
\begin{align*}
D =\begin{pmatrix}
  a & 0 & b \\
  0 & 0 & 0 \\
  b^* & 0 & c
\end{pmatrix}\begin{array}{c}
  \ese\ominus\ene \\
  \ene \\
  \ese^\bot
\end{array}&&\mbox{and}&&
D_k =\begin{pmatrix}
  a+\frac{1}{k}I & 0 & b \\
  0 & \frac{1}{k}I & 0 \\
  b^* & 0 & c+ \frac{1}{k}I
\end{pmatrix}\begin{array}{c}
  \ese\ominus\ene \\
  \ene \\
  \ese^\bot
  \end{array} ,
\end{align*}
where $a$, and therefore $\displaystyle a+\frac{1}{k}I$ are
invertible. Hence, by by Theorem \ref{PA}, 

\begin{align*}
\pas{D_k}{\ese} =\begin{pmatrix}
  I & 0 & (a+\frac{1}{k}I)^{-1}b \\
  0 & I & 0 \\
  0 & 0 & 0
\end{pmatrix}&&\mbox{and}&&\pas{D}{\ese} =
\begin{pmatrix}
  I & 0 & a^{-1}b \\
  0 & I & 0 \\
  0 & 0 & 0
\end{pmatrix}.
\end{align*}
So we obtain
\[
\|\pas{D}{\ese}\|=\lim_{k\rightarrow
\infty}\|\pas{D_k}{\ese}\|\leq \sup_{D'\in\edemas}\|\pas{D'}{\ese}\|
\]
which proves one inequality. The other inequality is a
consequence of equation (\ref{stewart}).\QED

\medskip

\begin{cor}\label{casofinito}
Let $\ese$ be a subspace of $\mathbb{C}^n$. Then
\begin{equation}\label{vuelta de stewart}
\sup_{D\in\edemas_n}\|\pas{D}{\ese}\|=\max
\Big{\{}\|\pas{Q}{\ese}\| : Q\in\edepn:\,R(Q)\bp\ese \Big{\}} .
\end{equation}
\end{cor}

\begin{rem}
This Corollary and the equation (\ref{StewartOLeary}) are connected
in the following sense. Observe that, by Theorem \ref{pqese menor que pdese},
the maximum can be taken either over all projections in position P' with $\ese$ or
over all the diagonal projections. So, if $I\subseteq\{1,\ldots,n\}$ and
$Q_I$ is the orthogonal projection onto the diagonal subspace spanned by
$\{e_i:\;i\in I\}$, then
\[
\sup_{D\in\edemas_n}\|\pas{D}{\ese}\|=\max
\Big{\{}\|\pas{Q_I}{\ese}\| : I\subseteq\{1,\ldots,n\}\Big{\}}
\]
Given a fixed $I\subseteq\{1,\ldots,n\}$, by Proposition
\ref{formula de ptak} it is easy to see that for a given diagonal
projection $Q$
\[
\|\pas{Q_I}{\ese}\|^{-2}=\sang{\ese}{N(Q_I)}=\min\{\pint{Q_I
x,\,x}:x\in \ese\ominus(\ese\cap N(Q_I))\;\mbox{and}\;\|x\|=1\}
\]
Since $(\ese\cap N(Q_I))\oplus \ese^\bot$ is the null space of
$P_\ese Q_I P_\ese$, the previous equality can be rewritten as
\[
\|\pas{Q_I}{\ese}\|^{-2}=\min\{\la\in \spec{P_\ese Q_I
P_\ese}:\;\la\neq 0\}=\min\{\la\in \spec{(Q_I P_\ese
Q_I}:\;\la\neq 0\}
\]
Consequently, if   $U$ is a matrix whose  columns form  an
orthonormal basis of $\ese$, then we get $Q_I P_\ese Q_I=Q_IUU^*Q_I$.
Therefore
\[
m_I^2=\min\{\la\in \spec{Q_IP_\ese Q_I}:\;\la\neq 0\}
\]
which shows that equations (\ref{StewartOLeary}) and (\ref{vuelta
de stewart}) are equivalent. In particular, note that  if $A$ is a
full rank matrix whose range is $\ese$ and the subspaces $R(Q_I)$
and $\ese$ are in position P', then
\[
\|\pas{Q_I}{\ese}\|=\|A(Q_IA)^{-1}Q_I\|=m_I^{-1}.
\]
\end{rem}

\bigskip

Next, we consider projections with complex weights. These
projections are studied by Bobrovnikova and Vavasis in
\cite{[vavasis2]}, who define, for each positive real number
$\mu$, the sets
\begin{align*}
C_\mu&=\{z\in\mathbb{C}:\;|\Pim z|\leq \mu \Preal
z\;\mbox{and}\;z\neq 0\} \  \hbox{ and } \   \ede_\mu
=\{D:\;D\in\ede\;\mbox{with entries in $C_\mu$}\},
\end{align*}
and they prove that
\[
\chinchulin{A}{\mu}=\sup_{D\in\ede_\mu}\|A(A^*DA)^{-1}A^*D\|<\infty.
\]

Let $Z$ be an $m\times (m-n)$ matrix, such that its columns form a
basis of $R(A)^\bot$. If the weights are positive, i.e. if $\mu=0$,
then Gonzaga and Lara \cite{[clovislara]} prove that
$\chinchulin{A}{0}=\chinchulin{Z}{0}$. In the next Proposition, we
generalize this result for $\mu\neq 0$.

\bigskip
\begin{pro}
Let $A$ be an $m\times n$ matrix of full rank, and let $Z$ be an
$m\times (m-n)$ matrix whose columns form a basis of $R(A)^\bot$.
Then,
\[
\chinchulin{A}{\mu}=\chinchulin{Z}{\mu}\hspace{1cm}\forall \mu>0
\]
\end{pro}
\dem Fix $D\in\ede_\mu$. On one hand, $A(A^*DA)^{-1}A^*D$ is
idempotent and
\begin{align*}
R(A(A^*DA)^{-1}A^*D)&=R(A)=:R\\
N(A(A^*DA)^{-1}A^*D)&=N(A^*D)=R(D^*A)^\bot=\Big[D^*(R(A))\Big]^\bot=:N\\
\intertext{On the other hand, $Z(Z^*D^{-1}Z)^{-1}Z^*D^{-1}$ is
also idempotent and}
R(Z(Z^*D^{-1}Z)^{-1}Z^*D^{-1})&=R(Z)=R(A)^\bot=R^\bot\\
N(Z(Z^*D^{-1}Z)^{-1}Z^*D^{-1})&=N(Z^*D^{-1})=R(D^{*-1}Z)^\bot
=\Big[D^{*-1}(R(A)^\bot)\Big]^\bot\\&=D(R(A))=N^\bot
\end{align*}

Using the fact that
$\|P_RP_N\|=\angf{R}{N}=\angf{R^\bot}{N^\bot}=\|P_{R^\bot}P_{N^\bot}\|$, 
and Ljance-Ptak's formula (Proposition \ref{formula de ptak}), we
obtain
\begin{align*}
\|A^*(ADA^*)^{-1}AD\|&=(1-\|P_{R}P_{N}\|^2)^{-1/2}\\
&=(1-\|P_{R^\bot}P_{N^\bot}\|^2)^{-1/2}=\|Z^*(ZD^{-1}Z^*)^{-1}ZD^{-1}\|
\end{align*}

Finally, since the map $D \rightarrow D^{-1}$
is a bijection of the set $\ede_\mu$, the result follows just by
taking supremum over all positive definite diagonal matrices.
\QED

\section{Compatibility of subspaces and orthonormal basis}

\subsection*{Definitions and main results}

Throughout this section, $\hil$ is a separable Hilbert space with
a fixed orthonormal basis $\bon$. Consider the abelian algebra $\ede$
of all operators which are diagonal with respect to  ${\mathcal B}$, i.e.
$C \in \ca $ belongs to $\ede$ if there exists a bounded sequence of
complex numbers $\{ c_n\}$ such that $C e_n = c_n e_n$  $(n \in \mathbb {N})$.
Denote by $\edemas$ the set of all positive invertible operators of $\ede$
and by $\edep$ the set of all projections of $\ede$.

Let us extend the definition of compatibility to the context of the diagonal algebra $\ede$.

\bigskip
\begin{fed}\label{definicion de compatibilidad}
A closed subspace $\ese$ of $\hil$  is
\textbf{compatible} with ${\mathcal B}$ (or \textbf{${\mathcal B}$-compatible}) if
\begin{equation}\label{condicion de compatibilidad con angulos}
  \sup_{Q\in\edep}\angf{\ese}{R(Q)}:= \angf{\ese}{\ede}<1 .
\end{equation}
In this case, we define $$ \con{\ese} {\ede } =
(1-\angf{\ese}{\ede}^2)\mrai = \left( \inf
_{Q\in\edep}\sang{\ese}{R(Q)}\right)\inv $$
\end{fed}

\begin{rems}\label{cuenta de los gaitas}
$\;$
\begin{description}
  \item[1.]Since $\angf{\ese}{\ete}=\angf{\ese^\bot}{\ete^{\bot}}$ for every pair of closed subspaces
  $\ese $ and $\ete$,  a subspace $\ese$ is ${\mathcal B}$-compatible if and only if
  $\ese^\bot$ is ${\mathcal B}$-compatible. Moreover, $\angf{\ese}{\ede} =
\angf{\ese\orto}{\ede}$.
  \item[2.] If  the dimension of the Hilbert space is finite,
  every subspace is compatible with every orthonormal basis.
\end{description}
\end{rems}

\bigskip
\noi
Using Remark \ref{equiva} we can give an alternative
characterization of compatibility.

\bigskip
\begin{teo}\label{condicion de compatibilidad con los pqs}
Let $\ese$ be a closed subspace of $\hil$. Then, $\ese$ is
${\mathcal B}$-compatible if and only if
\begin{equation}\label{equation de los pqs menor que infinito}
  \sup_{Q\in\edep}\|\pas{Q}{\ese}\|<\infty .
\end{equation}
Moreover, in this case, $$ \sup_{Q\in\edep}\|\pas{Q}{\ese}\| = K
[\ese , \ede] \quad . $$
\end{teo}
\dem Given a projection $Q$, by Remark \ref{equiva} we know that
  \[
  \|\pas{Q}{\ese}\|=(1-\angf{\ese}{N(Q)  }^2)\mrai.
  \]
If $\ese$ is compatible with ${\mathcal B}$, then
$\angf{\ese}{\ede}<1$, and therefore
\[
\sup_{Q\in\edep}\|\pas{Q}{\ese}\|\leq (1-\angf{\ese}{\ede}^2)\mrai
= \con{\ese }{ \ede } <\infty.
\]
Conversely, if (\ref{equation de los pqs menor que infinito})
holds, there exist $M>1$ such that
$\sup_{Q\in\edep}\|\pas{Q}{\ese}\|\leq M.$ Therefore,
\[
\sup_{Q\in\edep}\angf{\ese}{R(Q)} = (1-M^{-2})\mrai < 1 .
\]
\QED

\noi
The main result of this section is the following
theorem which is the natural generalization of Theorem \ref{bental teboulle con pases}
(or, more precisely, Corollary \ref{casofinito})
to the infinite dimensional setting.

\bigskip
\begin{teo}\label{teo posta posta}
Let $\ese$ be a closed subspace of $\hil$. Then, the following
statements are equivalent
\begin{description}
    \item[1.]$\ese$ is compatible with ${\mathcal B}$;
    \item[2.]$\displaystyle \sup_{D\in\edemas}\|\pas{D}{\ese}\| <
\infty$.
\end{description}
In this case, it holds $$ \quad \sup_{D\in\edemas}\|\pas{D}{\ese}\| =
\con{\ese }{ \ede }  = \sup_{Q\in\edep}\|\pas{Q}{\ese}\| .
$$
\end{teo}

\medskip
\noi The proof of this Theorem will be divided in several parts.
We start with a technical result.

\bigskip
\begin{lem}\label{lema tecnico para la implicacion facil}
Let $\ese$ be a closed subspace of $\hil$ and suppose that $
\sup_{D\in\edemas}\|\pas{D}{\ese}\|<\infty$. Then, for every
projection $Q\in\edep$ the pair $(Q, \ese )$ is compatible, so that the
oblique projection $\pas{Q}{\ese}$ is well defined for all $Q \in \edep$.
\end{lem}
\dem Let $Q\in\edep$  and consider the sequence
of invertible positive operators $\{D_k\}_{k\geq 1}$ defined by
$D_k=Q+\frac{1}{k}I.$ Since $D_k$ is invertible, the projection
$\pas{D_k}{\ese}$ is well defined. Moreover, by hypothesis we know
that $ \sup_{k\geq 1}\|\pas{D_k}{\ese}\|<\infty $. Therefore,  the
sequence $\{\pas{D_k}{\ese}\}$ has a limit point $P$ in the weak
operator topology (WOT) of $\ca$, because the unit ball of $\ca$ is
$WOT$-compact (see 5.1.3 of \cite{[KR]}). Moreover, if $\H$ is separable, the ball is
metrizable for the weak operator topology. Therefore, 
we can suppose that $\pas{D_k}{\ese}\convwot P $.

We shall prove that $P\in\Pas{Q}{\ese}$, that is, $P^2 = P$,
$R(P)=\ese$ and $QP=P^*Q$. The first two conditions follow from
the fact that, for every $k \in \mathbb {N}$,
\[
\pas{D_k}{\ese}=\begin{pmatrix}
  1 & x_k \\
  0 & 0
\end{pmatrix} \begin{array}{l}
    \ese \\
    \ese^\bot \
  \end{array}  , \ \ \hbox { so that } \ \
P=\begin{pmatrix}
  1 & x \\
  0 & 0
\end{pmatrix} \begin{array}{l}
    \ese \\
    \ese^\bot \
  \end{array} ,
\]
where $x$ is the $WOT$-limit of the sequence $x_k = P_\ese D_k
(1-P_\ese) $. On the other hand, for each $k \in \mathbb {N}$,
\[
D_k\pas{D_k}{\ese}=\pas{D_k}{\ese}^*D_k \ .
\]
An easy $ \frac{\varepsilon}{2}$ argument shows that
$D_k\pas{D_k}{\ese}\convwot Q P$, so, taking limit in the above
equality and using the fact that the involution is continuous in
the weak operator topology, we obtain $QP=P^*Q$.\QED

\medskip
\noi The next result, which can be proved in the same way 
as Proposition \ref{pqese menor que pdese},
by using Lemma \ref{lema tecnico para la implicacion facil},
provides the easier inequality in Theorem \ref{teo posta posta}:

\bigskip
\begin{pro}\label{implicacion facil}
Let $\ese$ be a closed subspace of $\hil$ and suppose that $
\sup_{D\in\edemas}\|\pas{D}{\ese}\|<\infty$. Then, $\ese$ is
${\mathcal B}$-compatible and
\begin{equation}\label{equation 2 implicacion facil con los pds}
\con{\ese}{ \ede} = \sup_{Q\in\edep}\|\pas{Q}{\ese}\|\leq
\sup_{D\in\edemas}\|\pas{D}{\ese}\|
\end{equation}
\end{pro}


\medskip

\noi
The other inequality of Theorem \ref{teo posta posta} is more complicated,
so we first need to  prove a particular case of it.
\medskip
\noi
For each $n\in \mathbb{N}$, denote
$\H_n = \generado{e_1, \dots, e_n}$,  $Q_n$ the  orthogonal
projection  onto $\H_n$, and for any closed
subspace $\ese$, denote $\ese_n = \ese \cap \H_n$. Recall that, given two
orthogonal projections $P$ and $Q$, $P\wedge Q$ denotes
the orthogonal projection onto $R(P)\cap R(Q)$.

\bigskip
\begin{pro}\label{bental cuando ese es finito}
Let $\ese$ be a finite dimensional  subspace of $\hil$ such that,
for some  $n \in \mathbb{N}$, $\ese \inc \H_n$.  Then $\ese$ is ${\mathcal B}$-compatible.
Moreover, if $E\in\edepfin$ satisfies $P_\ese\leq E$, then
\[
\{\pas{D}{\ese}:D\in\edemas\}\subseteq \cc
\{\pas{Q}{\ese}:Q\in\edep , \;Q\leq E\; \mbox{and}\;\, R(Q) \bp
\ese\} .
\]
In particular,
\begin{align*}
\sup_{D\in\edemas}\|\pas{D}{\ese}\|&\leq \sup
\{\|\pas{Q}{\ese}\|:Q\in\edep , \;Q\leq E\; \mbox{and}\;\,R(Q) \bp
\ese \} = \con{\ese}{\ede} .
\end{align*}
\end{pro}
\dem Let $E\in\edepfin$ be such that $P_\ese\leq E$. Denote $\ete=R(E)$.
Given $D\in\ede$, $D\ge 0$, the subspace
$\ete$ induces a matrix decomposition of $D$, 
\[
D=\bm{cc} D_1 & 0\\0&D_2\em
\begin{array}{c}     \ete \\    \ete^\bot \  \end{array}.
\]
If the pair $(D, \ese ) $ is compatible, it is easy to see that
the pair $(D_1, \ese ) $ is compatible in $L(\ete )$ and, by
Proposition \ref{compri},
\begin{equation}\label{eq de pas dentro del pas}
  \pas{D}{\ese}=\begin{pmatrix}
    \pas{D_1}{\ese} & 0 \\
               0            & 0 \
  \end{pmatrix}\begin{array}{c}
    \ete \\
    \ete^\bot \
  \end{array},
\end{equation}
where $\pas{D_1}{\ese}$ is considered as an operator in $L(\ete)$.
Since $\dim \ete < \infty$,  we deduce that $\ese$ is ${\mathcal
B}$-compatible. The other statements follow from Theorems
\ref{bental teboulle con pases} and \ref{condicion de compatibilidad con los pqs}. \QED

\bigskip
\begin{lem}\label{union de los esen es densa en ese}
Let $\ese$ be a closed subspace of $\hil$ such that
\begin{equation}\label{condicion para el lema tecnico}
c \ := \ \sup \ \{  \angf  { \ese }{ \H_n } : n \in \mathbb{N} \}
<1 .
\end{equation}
Then
\[
\overline{\bigcup_{n=1}^{\infty}\ese_n}=\ese \ .
\]
\end{lem}
\dem The assertion of the Lemma is equivalent to
$$
P_{\ese}\wedge Q_n \convsoti P_\ese  \ .
$$
Let $\xi\in\hil$ be a unit vector and let $\eps>0$. Let $k \in
\mathbb{N}$ such that $\displaystyle c^{2k-1}\leq\frac{\eps}{2}$.
By Proposition \ref{aproximacion a la interseccion},  for every
$n\geq 1$ it holds that
\[
\left\|\left(P_{\ese}Q_n\right)^k-P_{\ese}\wedge Q_n\right\|\leq
\frac{\eps}{2}.
\]
On the other hand, since $Q_nP_\ese\convsot P_\ese$ and the
function $f(x)=x^k$ is SOT-continuous on bounded sets (see, for example,
2.3.2 of \cite{[Pe]}),
there exists $n_0\geq 1$ such that, for every $n\geq n_0$,
\[
\left\|\left[\left(Q_n P_\ese\right)^k-P_{\ese}\right]\; \xi
\right\|<\frac{\eps}{2} \ .
\]
Therefore, for every $n\geq n_0$,
\begin{align*}
\left\|\left(P_\ese-P_{\ese}\wedge Q_n\right)\; \xi \right\|&\leq
\left\|\left[P_\ese-\left(P_{\ese}Q_n\right)^k\right]\; \xi
\right\|+ \left\|\left(\left(P_{\ese}Q_n\right)^k-P_{\ese}\wedge
Q_n\right)\; \xi \right\|<\eps  \ .
\end{align*}
\QED

\medskip
\noi Observe that, using Proposition \ref{bental cuando
ese es finito} and Lemma \ref{union de los esen es densa en ese},
we get the following characterization of
finite dimensional subspaces which are ${\mathcal B}$-compatible:

\bigskip
\begin{cor}\label{pocos}
Let $\ese $ be a finite dimensional subspace of $\hil$. Then
$\ese$ is ${\mathcal B}$-compatible \sii \ there exists $n \in
\mathbb {N}$ such that $\ese \inc \hil_n$.
\end{cor}

\bigskip
\begin{lem}\label{el supremo en los  chiquitos es menor}
Let $\ese$ be a ${\mathcal B}$-compatible subspace and $\ese_n =
\ese \cap \hil_n$, $n \in \mathbb {N}$. Then
\[\con{\ese_n }{ \ede } =
\sup_{Q\in\edep}\|\pas{Q}{\ese_n}\|\leq
\sup_{Q\in\edep}\|\pas{Q}{\ese}\| =\con{\ese }{ \ede }
\]
for every $n\geq 1$.
\end{lem}
\dem Using Corollary \ref{bental cuando ese es finito} and
Theorem \ref{implicacion facil}, we get
\[
\con{\ese_n }{ \ede } =\sup \{\|\pas{Q}{\ese_n}\|:Q\in\edep , \;
Q\leq Q_n\; \mbox{and}\;\, R(Q) \bp \ese_n \}.
\]
Thus, it suffices to prove the inequality for every $Q\in\edep$ such
that $R(Q)\bp\ese_n$ and  $Q\le Q_n$.
For each such $Q$ consider $\widehat{Q}=Q+(1-Q_n)\in\edep$. Then
$N(\widehat{Q}) = N(Q)  \cap R(Q_n)$, and
\begin{align*}
\angf{N(Q)  }{\ese_n}&=\sup\{\,|\pint{\xi , \, \eta}|:\;\xi \in
N(Q)  ,\;\eta\in\ese_n \;\,\mbox{and } \,  \|\xi\|=\|\eta\|=1 \}\\
&=\sup\{\,|\pint{\xi , \, \eta}|:\;\xi \in N(Q)  \cap
R(Q_n),\;\eta\in\ese_n \;\,\mbox{and } \  \|\xi\|=\|\eta\|=1 \}\\
&\leq \sup\{\,|\pint{\xi , \, \eta}|:\;\xi \in N(Q)  \cap
R(Q_n),\;\eta\in\ese \;\,\mbox{and } \ \|\xi\|=\|\eta\|=1\}\\ &
=  \sup\{\,|\pint{\xi , \, \eta}|:\;\xi\in
N(\widehat{Q}),\;\eta\in\ese \;\,\mbox{and } \
\|\xi\|=\|\eta\|=1\}\\&= \angf{N( \widehat{Q})}{\ese} .
\end{align*}
Therefore, using Remark  \ref{equiva}, we obtain
\[
\|\pas{Q}{\ese_n}\| = \sang {N(Q) }{\ese_n} \inv \le \sang
{N(Q)}{\ese} \inv = \sang
{N(\widehat Q)}{\ese} \inv
=\|\pas{\widehat{Q}}{\ese}\|,
\]
which proves the desired inequality.\QED

\medskip
\begin{rem}\rm\label{cofinito}
The statement of Lemma \ref{el supremo en los  chiquitos es menor}
can be rewritten as
\[
\con{\ese_n }{ \ede } \le \left( \inf
_{E\in\edep}\sang{\ese}{R(E)}\right)\inv \ .
\]
Actually, we  proved that it suffices to take the infimum  over the
projections $E \in \edepfin$. Indeed, it is enough to consider
$E = 1- \widehat Q$, where $\widehat Q$ are the projections which
appear in the proof of Lemma \ref{el supremo en los  chiquitos es menor}.
\end{rem}

\bigskip
\noi
We can now complete the proof of
Theorem \ref{teo posta posta}.

\bigskip
\begin{pro}\label{implicacion dificil}
Let $\ese$ be a ${\mathcal B}$-compatible subspace of $\hil$. Then
\begin{equation}\label{equation de la implicacion dificil}
\sup_{D\in\edemas}\|\pas{D}{\ese}\| \le \con{\ese }{ \ede } .
\end{equation}
\end{pro}
\dem Fix $D\in\edemas$. Recall that  $\|\cdot\|_D$ denotes the norm defined by
$\xi \rightarrow \|D\rai \xi\|$. Since
$D$ is invertible,  $\|\cdot\|_D$ is equivalent to  $\|\cdot\|$; 
thus, the union of the subspaces $\ese_n$ is dense in $\ese$
under both norms $\|\cdot\|_D$ and $\|\cdot\|$. Since $\pas{D}{\ese}$ (resp.
$\pas{D}{\ese_n}$ ) is the $D$-orthogonal projection onto the
subspace $\ese$ (resp. $\ese_n$),  then for every unit
vector $\xi\in\hil$
\[
\|(\pas{D}{\ese_n}-\pas{D}{\ese})\; \xi \|_D\conv 0.
\]
Using again the equivalence between the norms $\|\cdot\|_D$
and $\|\cdot\|$, we get
\[
\|(\pas{D}{\ese_n}-\pas{D}{\ese})\; \xi \|\conv 0.
\]
On the other hand, using Lemma \ref{el supremo en los  chiquitos
es menor} and Proposition \ref{bental cuando ese es finito},
it holds that, for each $n\in \mathbb {N}$,
\[
\|\pas{D}{\ese_n}\; \xi \|\leq \|\pas{D}{\ese_n}\|\leq \con{\ese_n
}{ \ede } \le \con{\ese }{ \ede } .
\sup_{Q\in\edep}\|\pas{Q}{\ese}\|.
\]
Thus,
\[
\|\pas{D}{\ese}\; \xi
\|=\lim_{n\rightarrow\infty}\|\pas{D}{\ese_n}\; \xi \|\leq
\con{\ese }{ \ede },
\]
which completes the proof. \QED

\subsection*{Alternative characterizations of compatibility}

In this section we add some characterizations of compatibility
which involve only finite dimensional diagonal subspaces.

Let us begin with a Proposition whose proof is connected
with the proof of Theorem \ref{teo posta posta}. We use the notations
$\H_n $,  $Q_n$ and $\ese_n$ (for a closed subspace $\ese$)
as before.

\bigskip
\begin{pro}\label{equi1}
Let $\ese$ be a closed subspace of $\hil$. Then, the following
statements are equivalent
\ben
    \item $\ese$ is ${\mathcal B}$-compatible;
    \medskip
    \item
    \ben
    \item [a.]    $\displaystyle \bigcup_{n=1}^{\infty}\ese_n$ is dense in $\ese$;
    \medskip
    \item [b.] there exists $M > 0$ such that  $\con{\ese_n }{ \ede } =
\sup_{Q\in\edep}\|\pas{Q}{\ese_n}\|\leq M$ for every $n \in
\mathbb{N}$. \een \een
\end{pro}
\dem
\begin{description}
    \item[$2\Rightarrow 1.$] It is a consequence of lemmas \ref{union de los esen es densa en ese}
    and \ref{el supremo en los  chiquitos es menor}.
    \item[$1\Rightarrow 2.$] Following the same argument as in Proposition \ref{implicacion dificil},
    we obtain that
    \[
        \sup_{D\in\edemas}\|\pas{D}{\ese}\| <\infty,
    \]
    and, by Proposition \ref{implicacion facil},
    $\ese$ is ${\mathcal B}$-compatible. \QED
\end{description}

\medskip
\noi Given $T\in \ca$, its reduced minimum modulus (see, e.g.,  \cite{[gordo]}), is defined by

\[
   \gamma(T) = \inf \{ \|T \xi \| : \xi \in N( T )\orto , \|\xi \| =
   1 \}.
\]
It is easy to see that $\gamma (T)>0$ if and only if $R(T)$ is
closed. By Proposition \ref{producto con rango cerrado}, if $A, B \in \ca$ have closed ranges, then
$$ \gamma (AB) >0 \iff c[N(A)   , R(B) ] <1 . $$
The following proposition describes a useful relation between
angles and the reduced minimum modulus of an operator.

\bigskip
\begin{pro}\label{gama y c}
Let $T \in \ca$ and let $P\in \ca$ be a projection with $R(P) =
\ese$. Suppose that $\gamma(T) >0$. Then
\begin{equation}\label{los gamas}
\gamma (T) (1-\angf{N(T)  }{ \ese }^2\ )\rai  \le \gamma (TP) \le
\|T\| (1-\angf{N(T)  }{ \ese }^2 \ )\rai .
\end{equation}
\end{pro}
\dem Note that $c[N(T)  , \ese] = c_0 [N(T)  , \ese \ominus (N(T)
\cap \ese )] =\| P_{N(T) } P_{\ese \ominus (N(T)  \cap \ese )}\|$,
by Proposition \ref{propiedades elementales de los angulos}. On the
other hand, $$ N(T) P = P\inv (N(T) )= P\inv (N(T)  \cap \ese )=
N( P) \oplus (N(T)  \cap \ese ) , $$ so that $$ (N(T) P ) \orto =
\ese \ominus (N(T)  \cap \ese ) \inc \ese =R(P) . $$ If $\xi \in
(N(T) P) \orto$, $\|T P \xi \| = \| T \xi \| = \| T ( P_{N(T)
\orto } \xi ) \|.$ Therefore, for every $\xi \in N(TP) \orto$, $$
\gamma (T) \|P_{N(T) \orto } \xi \| \le \|T P \xi \| \le \|T\|
\|P_{N(T) \orto } \xi \| . $$ Now, if $\| \xi \| = 1$, then $$
\|P_{N(T) \orto } \xi \|^2 = 1 - \| P_{N(T) } \xi \| ^2 \ge 1 - \|
P_{N(T) } P_{\ese \ominus (N(T)  \cap \ese )}\|^2 =1-c[N(T)  ,
\ese ]^2 , $$ since $\xi \in \ese \ominus (N(T)  \cap \ese )$.
This shows that $\gamma (T) (1-c[N(T)  , \ese ]\ )\rai  \le \gamma
(TP)$.

In order to prove the second inequality, consider a sequence $\xi_n$ of unit
vectors in $\ese \ominus (N(T)  \cap \ese )$ = $(N(T) P )\orto$
such that $\|P_{N(T) } \xi _n \| \to \| P_{N(T) } P_{\ese \ominus
(N(T)  \cap \ese )}\| = c[N(T)  , \ese ]$; then,
$$ \gamma (TP)^2 \le \|T P \xi_n \|^2  \le \|T\|^2 \|P_{N(T) \orto
} \xi_n \|^2  = \|T\|^2 (1 - \| P_{N(T) } \xi_n \| ^2) \to \|T\|^2
(1-c[N(T)  , \ese ]^2) . $$ \QED

\bigskip
\noi
Recall the notations $\edepfin = \{ Q \in \edep : Q$ has finite rank$\}$ and
$\edepfins = \{ Q \in \edepfin : R(Q) \cap \ese = \{0\} \ \}$.

\bigskip
\begin{pro}\label{equi2}
Let $\ese$ be a closed subspace of $\hil$.
Then, the following conditions are equivalent:
\begin{description}
\bigskip
    \item[1.] $\ese$ is ${\mathcal B}$-compatible, that is
    $\displaystyle\sup_{Q\in\edep}\angf{\ese}{R(Q)}<1 ;$
    \item[2.] $\displaystyle\sup_{Q\in\edepfin } \angf{\ese}{R(Q)}<1 ;$
    \item[3.]$\displaystyle \sup_{Q\in\edepfins }\angf{\ese}{R(Q)}<1 .$
\end{description}
\end{pro}
\dem It is clear that $1\Rightarrow 2 \Rightarrow 3$. In order to
prove that $3\Rightarrow 2$, let  $T\in \op$ such that
$N(T)=\ese$ and $\gamma(T)>0$. Given $Q\in \edepfin $ 
there exists $E\in \edepfins $ 
such that $E\leq Q$ and $R(TQ)=R(TE)$.
So, using equation (\ref{los gamas}), we obtain
\begin{align*}
\gamma (T)^2 (1-\angf{\ese }{R(E)}^2\ ) &\le \gamma
(TE)^2=\gamma(TET^*) \le \gamma(TQT^*) = \gamma(TQ)^2\\&\leq \|T\|^2
(1-\angf{\ese }{R(Q)}^2 \ ) .
\end{align*}
Therefore,
\[
0<\frac{\gamma (T)^2}{\|T\|^2}\ \
\inf_{Q\in\edepfins } (1-\angf{\ese }{R(E)}^2\
)\leq 1-\angf{\ese }{R(Q)}^2.
\]
Since we have chosen an arbitrary projection of
$\edepfin$ it holds that
\[
0<\inf_{Q\in\edepfin} (1-\angf{\ese }{R(Q)}^2).
\]
which is equivalent to
$\displaystyle\sup_{Q\in\edepfin }\angf{\ese}{R(Q)}<1 .$

Finally, if condition 2 holds,  we have that, in particular
\[
\sup \ \{  \angf  { \ese }{ \H_n } : n \in \mathbb{N} \}  <1,
\]
so $\displaystyle \overline{\bigcup_{n=1}^{\infty}\ese_n}=\ese$.
On the other hand, by Remark \ref{cofinito}, there exists $M\in
\mathbb{R}$ such that
\[
\con{\ese_n }{ \ede }<M \ , \hspace{1cm}\forall\, n\in\mathbb{N} .
\]
Therefore, by Proposition \ref{equi1}, $\ese$ is ${\mathcal
B}$-compatible. \QED

\section{An application to Riesz frames}\label{seccion de frames}

Recall that a sequence $\{ \xi_n \}$ of elements of $\H$ is called a $frame$
if there exist positive constants $A, B$ such that
\begin{equation}\label{frame}
A \|\xi \|^2  \le \sum |\api \xi , \xi _n \cpi |^2 \le B \|\xi
\|^2
\end{equation}
for all $\xi \in \H $.
The theory of frames, introduced by Duffin and Schafter \cite{[Duf]} in
their study of non harmonic Fourier series, has grown enormously 
after Daubechies, Grassmann and Meyer \cite{[DGM]} emphasized
their relevance in time-frequency analysis. The reader is referred
to the books by Young \cite{[young]} and Chistensen \cite{[Chrbook]}, and
the surveys by Heil and Walmut \cite{[HeiWal]} Casazza \cite{[Cas]}
and Christensen \cite{[Chris]} for
modern treatments of frame theory and applications.

It is well known that  $\{ \xi_n \}$ is a frame of $\H$ \iiff
there exists  an orthonormal basis $\{e_n\} _{n \in \mathbb {N}}$
of $\hil$ and an epimorphism (i.e. surjective) $T\in \ca$ such that
$\xi_n = Te_n $, $ n \in \mathbb {N}$. If the basis $\{e_n\}_{n \in \mathbb {N}}$
is fixed, the operator $T$ is called the $analysis$
operator and $T^*$, given by $T^*\xi = \sum_n \api \xi , \xi_n
\cpi e_n$, is called the $synthesis$ operator of the frame. The
positive inversible operator $S = TT^*$ (given by $S\xi =\sum_n
\api \xi , \xi_n \cpi \xi_n$) is called the $frame$ operator.
In this case, the optimal constants for equation (\ref{frame}) are
$B = \|S\| = \|T\|^2 $ and $A = \| S\inv\| \inv = \gamma (T)^2$.

The frame $\{ \xi_n\}$ is called a $Riesz$ $frame$ (see
\cite{[Chri]}) if there exists $C>0$ such that, for every $J \inc \mathbb {N}$,
the sequence $\{ \xi_n \} _{n \in J}$ is a frame (with
constants $A_J$ and $B_J$) of the space $\hil_J = \overline {\generado
{\xi_n : n \in J}}$  and $A_J \ge C$.

Consider $P_J = P_{\hil_J} \in \edep$. It is easy to see that $\{ \xi_n\}$
is a Riesz frame \sii \ there exists $c>0$ such that  $c \le
\gamma (TP_J )$ for  every $J \inc \mathbb {N}$. We prove now that
this condition is equivalent to the fact that $N(T) $ is
compatible with the basis $\{e_n\}_{n \in \mathbb {N}}$:


\bigskip
\begin{teo} Given an orthonormal basis $(e_n)_{n \in \mathbb {N}}$ of
$\hil$ and an epimorphism $T \in \ca$, then $(Te_n
)_{n \in \mathbb {N}}$ is a Riesz frame \sii \ $N(T) $ is
compatible with respect to the basis $(e_n)_{n \in \mathbb {N}}$.
\end{teo}
\dem Fix $J\subseteq \mathbb{N}$.
Then $R(TP_J )$ is closed \sii \ $\angf{N(T) } {\hil_J} <1$ (that
is, $\gamma(TP_J ) \neq 0$) and, in this case, $T|_{\hil_J} :
\hil_J \to R(TP_J )$ defines a frame with constants $A_J = \gamma
(TP_J )^2 $ and $B_J = \|TP_J\|^2$. Now, using Proposition
\ref{gama y c}, the statement becomes clear, because the frame
defined by $T$ is a Riesz frame \sii \ $ \inf _{J \inc  \mathbb
{N}} \gamma (TP_J) \ > \ 0$, which is equivalent to  $  \sup_{J \inc  \mathbb {N}}
\angf{N(T) }{\hil_J}  = \angf{N(T) }{\ede} < 1 . $ \QED

\bigskip
\begin{cor}\label{pocos riezses}
Given an  orthonormal basis $\{e_n\}_{n \in \mathbb {N}}$ of $\hil$
and an epimorphism $T \in \ca$ such that $N(T)$ has finite
dimension, then $\{Te_n \}_{n \in \mathbb {N}}$ is a Riesz frame
\sii \
 there exists $n \in \mathbb {N}$ such that $N(T) \inc \generado{e_1, \dots, e_n}$.
\end{cor}

\bigskip
\begin{rem}\rm Further applications of the relationship between Riesz frames and
compatible subspaces will be developed in a forthcoming paper.
\end{rem}



\end{document}